\documentclass[12pt]{amsart}

\usepackage{amsthm, amssymb, amsmath}

\theoremstyle{plain}
\numberwithin{equation}{section}
\newtheorem{theorem}[equation]{Theorem}
\newtheorem{prop}[equation]{Proposition}
\newtheorem{lemma}[equation]{Lemma}

\theoremstyle{definition}
\newtheorem{definition}[equation]{Definition}

\title{A Schwarz Lemma on the Polydisk}

\author{Greg Knese \\
        Washington University\\
        St. Louis, MO 63130}

\date{\today}

\begin{document}

\begin{abstract}
  We prove a generalization of the infinitesimal portion of the
  classical Schwarz lemma for functions from the polydisk to the disk.
  In particular, we describe the functions which play the role of
  automorphisms of the disk in this context--they turn out to be
  rational inner functions in the Schur-Agler class of the polydisk.
  In addition, some sufficient conditions are given for a function to
  be of this type.
\end{abstract}

\maketitle

\begin{section}{Introduction}
  The Schwarz Lemma on the unit disk $\mathbb{D} \subset \mathbb{C}$
  says that if $f: \mathbb{D} \to \mathbb{D}$ is a holomorphic
  function with $f(0)=0$, then $|f(z)|\leq |z|$ for all $z \in
  \mathbb{D}$ and $|f'(0)| \leq 1$.  Furthermore, if equality ever
  occurs in either of the two preceding inequalities, then $f$ is of
  the form $f(z)=\lambda z$ where $|\lambda| =1$.  Needless to say,
  this result has been generalized to many different situations.
  
  Here we would like to focus on the infinitesimal part of the Schwarz
  Lemma, which can be stated invariantly (without reference to the
  origin) as $$|f'(z)| \leq \frac{1-|f(z)|^2}{1-|z|^2}$$
  provided $f$
  is a holomorphic function from the disk to the disk.  Equality at a
  single point implies equality at every point and that $f$ is an
  automorphism of the disk:
  $$
  f(z) = \lambda \frac{z-a}{1-\bar{a} z} $$
  for some unimodular
  constant $\lambda$ and some point $a$ in the disk.  For holomorphic
  functions from the polydisk $\mathbb{D}^n=\mathbb{D}\times \cdots
  \times \mathbb{D}$ to the disk, we have the following well-known
  simple generalization (see \cite{RUD} page 179).

\begin{prop} If $f: \mathbb{D}^n \to \mathbb{D}$ is holomorphic then
\begin{equation}
  \sum_{i=1}^{n} (1-|z_i|^2) \left|\frac{\partial f}{\partial
      z_i}\right| \leq 1-|f(z)|^2
\label{sch}
\end{equation}
for any $z=(z_1, z_2, \dots, z_n) \in \mathbb{D}^n$.  
\end{prop}

In light of the situation for $n=1$, we feel compelled to ask, for
which functions $f:\mathbb{D}^n \to \mathbb{D}$ does equality in
(\ref{sch}) hold at every point?  The following theorem answers this
question.

\begin{theorem} Let $f: \mathbb{D}^n \to \mathbb{D}$ be a holomorphic
  function which depends on all $n$ variables (i.e. none of $f$'s
  partial derivatives is identically zero).
Then, $f$ satisfies 
\begin{equation}
  \sum_{i=1}^{n} (1-|z_i|^2) \left|\frac{\partial f}{\partial
      z_i}\right| = 1-|f(z)|^2
\label{equality}
\end{equation}
for every $z \in \mathbb{D}^n$ if and only if there is an $(n+1)\times
(n+1)$ \emph{symmetric} unitary $U$ (i.e. $U^*=U^{-1}$ and $U=U^t$)
where
\begin{equation*}
\begin{array}{cc} & \hspace{.07in} \begin{array}{cc} \mathbb{C} & \mathbb{C}^n
  \end{array} \\
U=\begin{array}{l} \mathbb{C} \\ \mathbb{C}^n \end{array} \hspace{-.2in}&
\left(\begin{array}{cc} A & B \\ C & D \end{array}\right) \end{array}
\end{equation*}
such that
\begin{equation}
f(z) = A + B E_z(I-D E_z)^{-1}C
\label{form}
\end{equation}
where $E_z$ is the $n\times n$ diagonal matrix with diagonal entries
$z_1, z_2, \dots, z_n$.  Such functions are rational and inner.
\label{thm1}
\end{theorem}

Interestingly, while the above extremal problem (\ref{equality}) is
posed entirely within the holomorphic functions from the polydisk to
the disk, its solution lies in the unit ball of the Schur-Agler class
of the polydisk.  Though the above functions are
generalizations of automorphisms of the disk, they are obviously not
automorphisms in any sense themselves.  Instead, the proper way to
think of them is as functions with the property that at every point of
the polydisk there is a holomorphic disk passing through that point on
which the restricted function is an automorphism of the disk.  Put
more simply, at every point of the polydisk the function is
``extremal'' in some direction.  Let us give a name to these
functions.

\begin{definition} A holomorphic function $f:\mathbb{D}^n \to
  \mathbb{D}$ is said to be \emph{everywhere infinitesimally
    extremal}, or \emph{extremal} for short, if $f$ satisfies
  (\ref{equality}) at every point of $\mathbb{D}^n$.
\end{definition}

The above theorem is proved using the technique commonly referred to
as the ``lurking isometry'' technique along with the ``polarization''
theorem for holomorphic functions, which (loosely stated) says that if
$F(z,\bar{z})\equiv 0$, then $F(z,w)\equiv 0$.

The next obvious question is on how big of a set does (\ref{equality})
need to hold in order for (\ref{equality}) to hold on the entire
polydisk?  A first guess might be that equality at one point implies
equality at every point.  The simple example $f(z_1,z_2)=(z_1+z_2)/2$ on the
bidisk proves that this is not the case.  Instead, we offer the
following two theorems as a substitute.

First, we make the following definition.

\begin{definition} Given a holomorphic function $f:\mathbb{D}^n \to
  \mathbb{D}$, we define the \emph{extreme set} $X(f) \subset
  \mathbb{D}^n$ for $f$ to be the set of points for which $f$
  satisfies equality in (\ref{equality}), i.e.
  $$
  X(f) := \left\{ z \in \mathbb{D}^n : 1-|f(z)|^2 = \sum_{i=1}^n
    (1-|z_i|^2) \left|\frac{\partial f}{\partial z_i}\right| \right\}
  $$
\end{definition} 

\begin{theorem} Let $f:\mathbb{D}^n \to \mathbb{D}$ be holomorphic.
  If $X(f)$ has nonempty interior, then $f$
  is everywhere infinitesimally extremal; i.e. $X(f)=\mathbb{D}^n$.
\label{thm2}
\end{theorem}

The above theorem is a direct result of the proof of theorem
(\ref{thm1}).  In order to state the next theorem, we need another
definition.

\begin{definition} A subset $D \subset \mathbb{D}^n$ is called a
  \emph{balanced disk} if $D$ is of the form
$$ D = \{(\phi_1(\zeta), \phi_2(\zeta), \dots, \phi_n (\zeta)) : \zeta
  \in \mathbb{D}\} $$
where $\phi_1, \phi_2, \dots, \phi_n$ are automorphisms of the disk.
\label{ddisk}
\end{definition}

The significance of balanced disks lies in the fact that holomorphic
functions on the polydisk are maximized along certain balanced disks,
just as the $\ell^1$ norm of a vector $v\in \mathbb{C}^n$ is attained
by $v\cdot w$ for some $w$ with unimodular entries.  The motivation
for the term ``balanced'' is that given two points $z=(z_1,\dots,
z_n), w=(w_1,\dots, w_n)$ in a balanced disk $D \subset \mathbb{D}^n$
we always have that $\rho (z_i, w_i)$ does not depend on $i$, where
$\rho$ is the pseudo-hyperbolic distance on the disk.  It turns out
that $X(f)$ is always a union of balanced disks, but as the following
theorem suggests, the structure of $X(f)$ is somewhat limited.

\begin{theorem} Let $f:\mathbb{D}^2 \to \mathbb{D}$ be holomorphic.
  If $X(f)$ contains two distinct intersecting balanced disks, $f$ is
  everywhere infinitesimally extremal.
\label{thm3}
\end{theorem}

In section 2, we will prove theorems (\ref{thm1}) and (\ref{thm2}); in
section 3, we will prove theorem (\ref{thm3}); and in section 4, we
will present some examples and curiosities.  Specifically, we present
an example (constructed using Nevanlinna-Pick interpolation on the
bidisk) which shows that it is possible for $X(f)$ to contain two
non-intersecting balanced disks but not be the whole bidisk.  Also, we
will explicitly describe the extremal functions on the bidisk which
map the origin to zero.  Except for the coordinate functions, the
extremal functions on the bidisk never extend continuously to the
closure of the bidisk. We do not know if this is the case in higher
dimensions.

\end{section} 

\begin{section}{Proofs of theorems (\ref{thm1}) and (\ref{thm2})}

\begin{proof}[Proof of theorems (\ref{thm1}) and (\ref{thm2}):]
Let $f:\mathbb{D}^n \to \mathbb{D}$ be holomorphic and suppose $f$
satisfies 
\begin{equation}
\sum_{i=1}^{n} (1-|z_i|^2) \left|\frac{\partial f}{\partial
      z_i}\right| = 1-|f(z)|^2
\label{equality2}
\end{equation}
on some open subset $V$ of the polydisk.  By assumption, $f$ depends
all $n$ variables.  So, we may find a point $\alpha$ in $V$ at which
none of $f$'s partial derivatives vanishes.  Shrinking $V$ if necessary,
we may assume $V$ is a ball on which none of $f$'s partial derivatives
vanishes.  For $i=1,2,\dots,n$, we can find a holomorphic square root
of $\frac{\partial f}{\partial z_i}$ on $V$; call each square root
$g_i$.  Then, (\ref{equality2}) can be rewritten as
\begin{equation*}
\sum_{i=1}^{n} (1-|z_i|^2) |g_i(z)|^2 = 1-|f(z)|^2
\end{equation*}
for each $z \in V$.  The above equation is much nicer because it
allows us to polarize.  Let $V^{*}$ denote the set $\{\bar{z}: z\in V\}$
and define a holomorphic function $F:V \times V^{*}\to \mathbb{C}$ by
\begin{equation*}
F(z,w) = (1-f(z)\bar{f}(\bar{w}))-\left(\sum_{i=1}^{n} (1-z_i w_i)
g_i(z) \bar{g}_i (\bar{w})\right)
\end{equation*}
Since $F(z, \bar{z})=0$ for all $z\in V$, we conclude that $F(z,w)=0$
for all $z\in V$ and $w \in V^{*}$. Replacing $w$ with $\bar{w}$, we
have that for any $z, w\in V$
\begin{equation}
(1-f(z)\bar{f}(w))=\sum_{i=1}^{n} (1-z_i \bar{w}_i)
g_i(z) \bar{g}_i (w)
\label{polarized}
\end{equation}
and with a little rearranging we have
\begin{equation*}
1+\sum_{i=1}^{n} z_ig_i(z) \bar{w}_i \bar{g}_i(w) =
f(z)\bar{f}(w)+\sum_{i=1}^{n} g_i (z)\bar{g}_i(w)
\end{equation*}

The map which for each $z \in V$ sends
\begin{equation}
\left(\begin{array}{c} 1 \\ z_1 g_1 (z) \\ \vdots \\ z_n g_n (z)
  \end{array} \right)
\mapsto \left( \begin{array}{c} f(z) \\ g_1 (z) \\ \vdots \\ g_n (z)
  \end{array} \right)
\label{map}
\end{equation}
defines a unitary $U$ on the span of elements of $\mathbb{C}^{n+1}$ of the
above form on the left into the span of elements of $\mathbb{C}^{n+1}$
of the above form on the right.  Extending $U$ if necessary to all of
$\mathbb{C}^{n+1}$ we get an $(n+1)\times (n+1)$ unitary, which we
write in the following block form, in order to derive what this
says about $f$:
\begin{equation*}
\begin{array}{cc} & \hspace{.07in} \begin{array}{cc} \mathbb{C} & \mathbb{C}^n
  \end{array} \\
U=\begin{array}{l} \mathbb{C} \\ \mathbb{C}^n \end{array} \hspace{-.2in}&
\left(\begin{array}{cc} A & B \\ C & D \end{array}\right). \end{array}
\end{equation*}
Let $E_z$ be the $n\times n$ diagonal matrix with diagonal entries
$z_1, z_2, \dots, z_n$ and let $G(z)$ be the column vector with
entries $g_1 (z), g_2 (z), \dots, g_n (z)$.  Then, (\ref{map}) says
\begin{eqnarray*}
A + BE_z G(z) & = & f(z) \\
C + DE_z G(z) & = & G(z)  
\end{eqnarray*}
which, first solving for $G(z)$ and then for $f(z)$, yields
\begin{eqnarray}
G(z) & = & (I-DE_z)^{-1} C \label{G} \\
f(z) & = & A+BE_z (I-DE_z)^{-1} C. \label{f} 
\end{eqnarray}

for any $z \in V$.  Since $V$ is an open set, equation (\ref{f}) must
therefore hold on the entire polydisk.  

It still must be shown that $U$ is symmetric; i.e. $U^t = U$, which
admittedly is a strange property for a unitary to have.  Our first
step will be to take the partial derivatives of $f$ from (\ref{f}) and
set those equal to the squares of the entries of (\ref{G}) in order to
deduce that $B^t = C$ and $(D^t-D)E_z G(z) =0$ for all $z \in
\mathbb{D}^n$.  Our second step will be to deduce that the set $\{ E_z
G(z): z \in \mathbb{D}^n \}$ spans $\mathbb{C}^n$, from which we may
conclude $D=D^t$ and more importantly that $U = U^t$.

Let $e_1, e_2, \dots, e_n$ denote the standard basis (column) vectors
of $\mathbb{C}^n$.  Then, by (\ref{f})
\begin{eqnarray}
\frac{\partial f}{\partial z_i} &=& B(I-E_z D)^{-1} E_{e_i}
(I-DE_z)^{-1} C \nonumber \\
 &=& (B(I-E_z D)^{-1} e_i) (e_i^t (I-DE_z)^{-1}C) \nonumber \\
&=& (e_i^t (I-D^t E_z)^{-1} B^t) (e_i^t (I-DE_z)^{-1} C). \label{partial1}
\end{eqnarray}
On the other hand, by (\ref{G})
\begin{equation}
\frac{\partial f}{\partial z_i} = (e_i^t (I-DE_z)^{-1} C)^2.
\label{partial2}
\end{equation}
Combining (\ref{partial1}) and (\ref{partial2}) we get
\begin{eqnarray}
e_i^t (I-DE_z)^{-1} C & = & e_i^t (I-D^t E_z)^{-1} B^t \text{ for all
  $i$, and so} \nonumber \\
(I-DE_z)^{-1}C & = & (I-D^t E_z)^{-1} B^t \label{fact1}
\end{eqnarray}
Setting $z = 0$ we see that $B^t = C$. Now, multiplying (\ref{fact1})
by $(I-D^t E_z)$ yields
\begin{equation*}
 C = (I-D^t E_z)(I-DE_z)^{-1} C = (I-DE_z)^{-1}C -D^t E_z
 (I-DE_z)^{-1} C
\end{equation*}
which implies
\begin{equation*}
D^t E_z (I-DE_z)^{-1} C = (I-DE_z)^{-1} C -C = DE_z (I-DE_z)^{-1} C
\end{equation*}
which finally yields
\begin{equation*}
(D^t-D)E_z G(z) = 0
\end{equation*}
for all $z \in \mathbb{D}^n$.

Next, we prove that as $z$ varies over $\mathbb{D}^n$, $E_z G(z)$
spans all of $\mathbb{C}^n$.  This would be extremely easy if none of
$G$'s entries were zero at the origin.  For then, we could set
$z=re_i$ and for $r$ small enough $E_z G(z)$ would be a nonzero
multiple of $e_i$.  Unfortunately, it might be the case that some of
$G$'s entries vanish at the origin.  We do know that none of $G$'s
entries vanish at the point $\alpha \in V$ (from way back at the
beginning of the proof).  So, now the strategy is to somehow shift our
attention to $\alpha=(\alpha_1, \alpha_2,\dots, \alpha_n)$ and try to
apply a similar argument.

Let $\phi_i(\zeta) = (\alpha_i - \zeta)/(1-\bar{\alpha}_i \zeta)$ for
$i=1, 2, \dots, n$.  Using the identity
$$ 1-z_i \bar{w}_i = \frac{(1-\phi_i (z_i)\bar{\phi}_i
  (w_i))(1-\bar{\alpha}_i z_i)(1-\alpha_i \bar{w}_i)}{1-|\alpha_i|^2}
$$
the equation (\ref{polarized}) can be rewritten as
\begin{equation*}
(1-f(z)\bar{f}(w))=\sum_{i=1}^{n} \frac{(1-\phi_i (z_i)\bar{\phi}_i
  (w_i))(1-\bar{\alpha}_i z_i)(1-\alpha_i \bar{w}_i)}{1-|\alpha_i|^2}
g_i(z) \bar{g}_i (w)
\end{equation*}
and we get a new $(n+1)\times (n+1)$ unitary $\tilde{U}$ which maps
\begin{equation}
\left(\begin{array}{c} 1 \\ \frac{\alpha_1-z_1}{\sqrt{1-|\alpha_1|^2}}
    g_1 (z) \\ 
\vdots \\ 
\frac{\alpha_n-z_n}{\sqrt{1-|\alpha_n|^2}} g_n (z)
  \end{array} \right)
\mapsto \left( \begin{array}{c} f(z) \\
    \frac{1-\bar{\alpha}_1z_1}{\sqrt{1-|\alpha_1|^2}} g_1 (z) \\
 \vdots \\ 
\frac{1-\bar{\alpha}_n z_n}{\sqrt{1-|\alpha_n|^2}}g_n (z)
  \end{array} \right).
\label{mapRewritten}
\end{equation}
This unitary is, in fact, uniquely determined because, looking at the
left hand side, and setting in turn $z=\alpha, \alpha+re_1, \dots,
\alpha+re_n$ for $r$ small enough, we see that the left hand side
spans all of $\mathbb{C}^{n+1}$ as $z$ varies over the polydisk since
$g_i(\alpha)\ne 0$ for $i=1, 2, \dots, n$.  Therefore, the vectors on
the right hand side of (\ref{mapRewritten}) span all of
$\mathbb{C}^{n+1}$.  This implies vectors of the form
$$
\left( \begin{array}{c} f(z) \\
    (1-\bar{\alpha}_1z_1) g_1 (z) \\
 \vdots \\ 
(1-\bar{\alpha}_n z_n) g_n (z)
  \end{array} \right)
$$
span $\mathbb{C}^{n+1}$.  But, a vector of the above form is the image
of
\begin{equation}
\left( \begin{array}{c} 1 \\
    z_1 g_1 (z) \\
    \vdots \\
    z_n g_n (z)
  \end{array} \right)
\label{1EZG}
\end{equation}
under the invertible matrix $U-T$, where $T$ is the diagonal matrix
with diagonal $0, \bar{\alpha}_1, \bar{\alpha}_2, \dots, \bar{\alpha}_n$.  Hence,
vectors of the form (\ref{1EZG}) span $\mathbb{C}^{n+1}$, and
consequently, vectors of the form $E_z G(z)$ span $\mathbb{C}^n$, as
desired.

This establishes the fact that $U=U^t$.

Conversely, suppose $f$ is of the form (\ref{f}) for some symmetric
unitary.  It should be noted that $f$ does indeed map into the disk,
because of the formula:
\begin{equation}
1-|f(z)|^2 = C^* (I-E_{z}^* D^*)^{-1} (I-E_{z}^* E_z) (I-DE_z)^{-1} C
\label{fformula}
\end{equation}
which holds regardless of whether or not $U$ is symmetric.  But,
assuming $U$ is symmetric, the above simplifies to
$$
1-|f(z)|^2 = \sum_{i=1}^n (1-|z_i|^2) |e_{i}^t (I-DE_z)^{-1} C|^2.
$$

Also, since $U$ is symmetric, the right hand side is equal to

$$
\sum_{i=1}^n (1-|z_i|^2)\left|\frac{\partial f}{\partial z_i}\right|
$$
by (\ref{partial1}), which was derived straight from (\ref{f}).

Finally, $f$'s of the form (\ref{f}) are rational of degree
at most $2n$ by Cramer's rule.  They are inner by (\ref{fformula})
because a rational function cannot vanish on a set of positive measure
of the boundary of $\mathbb{D}^n$. 

We have also established theorem (\ref{thm2}),
because we proved that if (\ref{equality2}) holds on an open set, then
$f$ is of the form (\ref{form}), which in turn implies that
(\ref{equality}) holds everywhere.  

In closing the proof, we remark that the condition that $f$ depends on
all $n$ variables is simply there to make the theorem statement nicer.
If $f$ doesn't depend on some variable, then $f$ can be treated as a
function on a lower dimensional polydisk.
\end{proof}
\end{section}

\begin{section}{Proof of theorem (\ref{thm3})}

For the rest of the paper we will restrict our attention to the case
$n=2$; we will write points of $\mathbb{D}^2$ as $(z,w)$, and we will
write partial derivatives of $f$ as $f_1$ and $f_2$.  We
recall the following definitions from the introduction restricted to
the case $n=2$.

\begin{definition} Let $f:\mathbb{D}^2 \to \mathbb{D}$ be
  holomorphic. Then, the \emph{extreme set} of $f$ is the subset of
  $\mathbb{D}^2$ defined by
\begin{equation}
X(f) := \left\{ (z,w) \in \mathbb{D}^2 : 1-|f|^2 
= (1-|z|^2) |f_1| +
(1-|w|^2) |f_2| \right\}
\label{extremeset}
\end{equation}
\end{definition}

\begin{definition} A subset $D \subset \mathbb{D}^2$ is called a
  \emph{balanced disk} if $D$ is of the form
$$ D = \{(\phi_1(\zeta), \phi_2(\zeta)) : \zeta
  \in \mathbb{D}\} $$
where $\phi_1$ and  $\phi_2$ are automorphisms of the disk.
\end{definition}

Observe that two balanced disks intersect at no points, one point, or
all points.

\begin{lemma} Let $f:\mathbb{D}^2\to \mathbb{D}$ be holomorphic.  
If either of $f$'s partial derivatives equals zero at a point of
  $X(f)$, then $f$ depends on only one variable and
  $X(f)=\mathbb{D}^2$.  Furthermore, if $f$ depends on both variables,
  then through each point $P$ of $X(f)$ there is a unique balanced disk
  $D$ contained in $X(f)$ and containing $P$ such that the restriction
  of $f$ to $D$ is an automorphism of the disk.  
\label{lemma1}
\end{lemma}

Here, $f$ restricted to a disk $D=(\phi_1, \phi_2)(\mathbb{D})$ (where
$\phi_1$ and $\phi_2$ are automorphisms of the unit disk) being an
automorphism of the unit disk just means $f\circ (\phi_1, \phi_2)$ is
an automorphism.

\begin{proof}
  It is a simple calculation to verify that if $P\in X(f)$, $\phi$ is
  an automorphism of the disk, and $\Phi$ is an automorphism of the
  bidisk, then $\Phi^{-1}(P) \in X(\phi \circ f \circ \Phi)$.
  Therefore, to prove the first part of the lemma, it is sufficient to
  prove $f(z,w)=w$ under the assumptions $(0,0)\in X(f)$, $f(0,0)=0$,
  $f_1(0,0)=0$ and $f_2(0,0)=1$ by applying appropriate automorphisms.
  For any $(z,w) \in \mathbb{D}^2$ with $|w| > |z|$ we see that the
  function from the disk to the disk defined by
  $g(\zeta)=f((z/w)\zeta, \zeta)$ satisfies $g'(0)=1$ and $g(0)=0$.
  By the classical Schwarz lemma, $g(\zeta)=\zeta$ for all $\zeta \in
  \mathbb{D}$ and in particular $f(z,w)=w$.  It follows that this
  holds on the entire bidisk.

By the first part of the lemma, if $f$ depends on both variables, then
neither of $f$'s partial derivatives can vanish at a point of $X(f)$.
To prove the second portion of the lemma, we can assume like before
that $P=(0,0)$, $f(0,0)=0$, $f_1(0,0)=a$, and $f_2(0,0)=b$ where $a$
and $b$ are positive real numbers with $a+b=1$.  A balanced disk
through the origin can be written in the form $\{(\zeta, \mu \zeta):
\zeta \in \mathbb{D}\}$ where $|\mu|=1$.  Now, $h(\zeta)=f(\zeta, \mu \zeta)$
is an automorphism of the disk precisely when $|a+\mu b|=1$ or
equivalently when $\mu=1$.  Therefore, $D=\{(\zeta,\zeta)\}$ is the unique balanced
disk containing $P$ on which $f$ is an automorphism of the disk.
Finally, $D \subset X(f)$ because if $h(\zeta)=f(\zeta,\zeta)=\zeta$, then
\[
1=|h'(\zeta)|=|f_1(\zeta,\zeta)+f_2(\zeta,\zeta)|\leq
|f_1(\zeta,\zeta)|+|f_2(\zeta,\zeta)| \leq 1
\]
which implies $(\zeta,\zeta)\in X(f)$ for all $\zeta \in \mathbb{D}$.
\end{proof}

In order to prove theorem (\ref{thm3}), we need to give a name to a special
class of balanced disks associated to $f$.

\begin{definition} A balanced disk $D$ is called an
  \emph{extreme} disk for $f$ if $f$ restricted
  to $D$ is an automorphism of the unit disk $\mathbb{D}$.
\end{definition}

A result of lemma (\ref{lemma1}) is the following: if $f:\mathbb{D}^2
 \to \mathbb{D}$ has two intersecting extreme disks, then
 $X(f)=\mathbb{D}^2$.  We are now ready to prove theorem (\ref{thm3}).  

\begin{proof}[Proof of Theorem (\ref{thm3}):]
  We may assume the intersection point is the origin and the two disks
  are $D_1 = \{(\zeta, \zeta): \zeta \in \mathbb{D}\}$ and $D_2 =
  \{(\zeta, \mu \zeta): \zeta \in \mathbb{D}\}$ where $|\mu|=1$.  We
  may also assume $D_1$ is an extreme disk, because $f$ has an extreme
  disk through the origin, and we might as well assume $D_1$ is it.
  Also, assuming $f$ depends on both variables, which we do, this
  implies $D_2$ is not an extreme disk.
  
  Now, let $P=(a, b)$ be some point of $X(f)$.  We can actually derive
  a formula for the extreme disk passing through $P$.  Let $\phi_c
  :\mathbb{D}\to \mathbb{D}$ be the automorphism of the disk for
  which $0 \mapsto c$, $c \mapsto 0$, and $\phi_c=\phi_c^{-1}$; i.e.
$$
\phi_c (\zeta) = \frac{c- \zeta}{1-\bar{c} \zeta}.
$$

Then, the derivative of $g(\zeta)=f(\phi_a(\zeta), \phi_b(\lambda
\zeta))$ at $0$ satisfies
\begin{eqnarray*}
  |g'(0)| &=& \left|f_1(P)(1-|a|^2)+\lambda f_2(P)(1-|b|^2)\right| \\
& = & |f_1(P)|(1-|a|^2)+|f_2(P)|(1-|b|^2) \\
& = & 1-|g(0)|^2
\end{eqnarray*}
exactly when $\lambda$ satisfies
\begin{equation}
\lambda = {\rm sgn}\left[ \frac{f_1(P)}{f_2(P)}\right]
\label{lambda}
\end{equation}
(which makes sense because neither partial derivatives can vanish at a
point of $X(f)$ in the bidisk without vanishing everywhere) in which
case $g$ is an automorphism of the disk and
$$
D(P) = \{ (\phi_a(\zeta), \phi_b (\lambda \zeta)): \zeta \in
\mathbb{D}\}
$$
is the extreme disk passing through $P$.
So, defining $\lambda (P)$ for each $P \in X(f)$ as in (\ref{lambda}),
we can define a map $F$ from $\mathbb{D}^2$ into $X(f)$ by
mapping
$$
(z, w) \mapsto (\phi_{z} (w), \phi_{\mu z} (\lambda(z,\mu
z) w).
$$

The map $F$ is injective for the following reasons.  If
$F(z,w)=F(z',w')$, then the extreme disk through $(z, \mu z)$
intersects the extreme disk through $(z', \mu z')$.  These two extreme
disks must therefore be the same, since we are assuming $f$ depends on
both variables.  This implies $z=z'$, because otherwise $D_2$ would
contain two points of the same extreme disk and since $D_2$ is
balanced this would force $D_2$ to be an extreme disk, which it is
not.  This immediately gives $w=w'$.  Since $F$ is continuous (or even
smooth), we see that $X(f)$ has nonempty interior and therefore
$X(f)=\mathbb{D}^2$.
\end{proof}

\end{section}

\begin{section}{Examples and Curiosities}

\begin{subsection}{Example: A Function with Two Non-Intersecting Extreme Disks}

Define $g:\mathbb{D}^2 \to \mathbb{D}$ by
\begin{eqnarray*}
g(z,w) & = &
\left\{ (z+w)/2 + (\sqrt{5}/6) (z^2-w^2) - (1/9) (z^2+7zw+w^2)\right. \\
 & &- \left. (\sqrt{5}/3) (z^2w-zw^2) - (1/2) (z^2w+zw^2) + z^2w^2 \right\}\\
 & \div & \left\{ 1 - (1/2)(z+w) + (\sqrt{5}/3) (z-w) - (1/9) (z^2+7zw+w^2) \right. \\
 & &- \left.  (\sqrt{5}/6) (z^2-w^2) + (1/2) (z^2w+zw^2)\right\}
\end{eqnarray*}

This function $g$ has two non-intersecting extreme disks since
\begin{equation*}
g(z,z) = z
\end{equation*}
\begin{equation*}
g\left( \frac{z-\sqrt{5}/3}{1-(\sqrt{5}/3)z} ,
  \frac{z+\sqrt{5}/3}{1+(\sqrt{5}/3)z} \right) =  z
\end{equation*}
and 

$$
\frac{z-\sqrt{5}/3}{1-(\sqrt{5}/3)z} \ne
\frac{z+\sqrt{5}/3}{1+(\sqrt{5}/3)z} \text{ for all $z \in \mathbb{D}$.}
$$

However, $g$ does not satisfy (\ref{equality}) everywhere.

The function $g$ was constructed using the Nevanlinna-Pick
interpolation theorem on the bidisk by looking for solutions to the
four point problem
\begin{eqnarray*}
(0,0) &\mapsto & 0 \\
(r, r) &\mapsto & r\\
(0,b) &\mapsto & \alpha \\
(a,0) &\mapsto & \beta 
\end{eqnarray*}
where $r<1$ and $|b|=|a|=|(\alpha-\beta)/(1-\bar{\alpha}\beta)|<1$.
(This is just a pair of extremal two-point problems.) For the
statement of the appropriate theorem see \cite{PICK} chapter 11.
\end{subsection}

\begin{subsection}{Everywhere Infinitesimally Extremal Functions on
    the Bidisk}
  
  Using theorem (\ref{thm1}) it is a relatively simple task to write
  down explicitly the everywhere infinitesimally extremal functions on
  the bidisk which map the origin to zero.  They are of the following
  form:
\[
f(z,w)=\mu \frac{az+bw-zw}{1-\bar{b}z-\bar{a}w}
\]
where $|\mu|=|a|+|b|=1$.  Notice that unless $a$ or $b$ equals zero $f$
never extends continuously to the closed bidisk.  Indeed, supposing
$a$ and $b$ are both positive real numbers and $\mu=1$ (after composing $f$ with
an automorphism), we see that
\[
\lim_{\theta\to 0} f(e^{i\theta}, e^{i t\theta})= 
\begin{cases}   1  & \text{if } t \ne -b/a \\
               -1  & \text{if } t = -b/a \end{cases}
\]

This immediately implies that on the bidisk the only extremal
functions that extend continuously to the closed bidisk depend on one
variable and are simply automorphisms of the disk (in one variable);
i.e. $f(z,w)=\phi(z)$ or $f(z,w)=\phi(w)$ for some automorphism of the
disk, $\phi$.

\end{subsection}

\end{section}

\end{document}